\documentclass[12pt]{article}%
\usepackage{amsfonts}
\usepackage{amsmath}
\usepackage{amssymb}
\usepackage{graphicx}%
\providecommand{\U}[1]{\protect\rule{.1in}{.1in}}
\newtheorem{theorem}{Theorem}

\newtheorem{example}[theorem]{Example}

\newtheorem{lemma}[theorem]{Lemma}

\newenvironment{proof}[1][Proof]{\noindent\textbf{#1.} }{\ \rule{0.5em}{0.5em}}
\begin{document}

\title{About the relatively uniform completion of a principal ideal}
\author{Elmiloud Chil\\Institut pr\'{e}paratoire aux \'{e}tudes d'ing\'{e}nieurs de Tunis \\2 Rue Jawaher lel Nahrou Montfleury 1008 TUNISIA\\Univesity of Tunis Tunisia\\Elmiloud.chil@ipeit.rnu.tn}
\date{}
\maketitle

\begin{abstract}
In [2] the author claims to provide a counterexample to a result in a recent
paper [1]. In this note, we prove that the details of his example is false and
this example is compatible with our result in [1] and so is not a countreexample.

\end{abstract}



A Riesz space $E$ is called \textit{Archimedean} if for each non zero $a\in E$
the set $\left\{  na,n=\pm1,\pm2,...\right\}  $ has no upper bound in $E$. In
order to avoid unnecessary repetition we will assume throughout the paper that
all Riesz spaces under consideration are Archimedean.

A Riesz space is called \textit{Dedekind complete} whenever every nonempty
subset that is bounded from above has a supremum. Let us say that a vector
subspace $G$ of a Riesz space $E$ is \textit{majorizing} $E$ whenever for each
$x\in E$ there exists some $y\in G$ with $x\leq y$.(or, equivalently whenever
for each $x\in E$ there exists some $y\in G$ with $y\leq x$)$.$ A Dedekind
complete Riesz space $L$ is said to be a \textit{Dedekind completion} of the
Riesz space $E$ whenever $E$ is Riesz isomorphic to an order dense majorizing
Riesz subspace of $L$ (which we identify with $E$). It is a classical result
that every Archimedean Riesz space $E$ has a Dedekind completion, which we
shall denoted by $E^{\delta}$.

The relatively uniform topology on Riesz spaces plays a key role in the
context of this work. Let us recall the definition of the \textit{relatively
uniform convergence}. Let $E$ be an Archimedean Riesz space and an element
$u\in E^{+}$. A sequence $(x_{n})_{n}$ of elements of $E$ \textit{converges
u-uniformly} to the element $x\in E$ whenever, for every $\epsilon>0$, there
exists a natural number $N_{\epsilon}$ such that $|x_{n}-x|\leq\epsilon u$
holds for all all $n\geq N_{\epsilon}$. In such a case we shall write
$x_{n}\rightarrow x(u)$. The element $u$ is called the \textit{regulator of
the convergence}. The sequence $(x_{n})_{n}$ \textit{converges relatively
uniformly} to $x\in E$, whenever $x_{n}\rightarrow x(u)$ for some $u\in E^{+}%
$. We shall write $x_{n}\rightarrow x$(r-u) if we do not want to specify the
regulator. \textit{Relatively uniform limits} are unique if and only if $E$ is
Archimedean. A nonempty subset $D$ of $E$ is said to be \textit{relatively
uniformly closed} if every relatively uniformly convergent sequence in $D$ has
its limit also in $D$. We emphasize that the regulator does need not to be an
element of $D$. The relatively uniformly closed subsets are the closed sets of
a topology in $E$, namely, the relatively uniform topology. The closure, with
respect to the relatively uniform topology of the Riesz space $E$ in its
Dedekind completion is a uniformly complete vector lattice, denoted by
$E^{ru}$ and referred to as the uniform completion of $E$. Thus $E$ is an
order dense majorizing Riesz subspace of $E^{ru}$.

\section{Comments}

Following a note published on arXiv [2] concerning the content of my research
findings [1], I consider it both important and necessary for the sake of
scientific integrity to respond to certain comments raised therein.

In [1] we prove the following result

\begin{lemma}
Let $E$ be an Archimedean Riesz space, $0\neq x\in E$, $E_{x}$ be the
principal ideal generated by $x$ in $E$ . Then
\[
(E^{ru})_{x}=(E_{x})^{ru}.
\]

\end{lemma}

Following we give the details that the author in [2] gave some false arguments.

\begin{example}
Consider the Archimedean Riesz space $E$ consisting of all linear piecewise
continuous functions on $[0,1]$ . Let $u\in E$ be defined by $u(t)=t$ for all
$t\in\lbrack0,1]$ . We have $E^{ru}=C[0,1]$ and then $(E^{ru})_{u}$ is the
principal ideal of $C[0,1]$ generated by $u$. We claim that
\[
(E_{u})^{ru}=(E^{ru})_{u}.
\]

\end{example}

\begin{proof}
First as it is noted in [1] we have $(E^{\delta})_{u}=(E_{u})^{\delta}.$ From
this equality the $(r.u)$ convergence in $E_{u}$ coincide with the one on
$(E^{\delta})_{u}$ which is the uniformly convergence on $[0,1].$ On the other
hand we remark that the inclusion $(E_{u})^{ru}\subset(E^{ru})_{u}$ is
obvious. For the converse inclusion, let $f\in(E^{ru})_{u}=\left\{  f\in
C[0,1]:\exists M\geq0,\left\vert f\right\vert \leq Mu\right\}  $ then
$f(0)=0.$ Now, by using a classical approximation Theorem there exists
$(p_{n})_{n}\in E$ which converges uniformly to $f$ on $[0,1]$ then
$p_{n}(0)\rightarrow f(0)=0$ so by replacing $p_{n}$ by $p_{n}-p_{n}(0)$ which
converges uniformly also to $f$ we can assume without loss of generality that
$p_{n}(0)=0$ which implies that $\underset{t\rightarrow0^{+}}{\lim}\frac
{p_{n}(t)}{t}=p_{n}^{\prime}(0^{+})$ the derived of $p_{n}$ on the right of
$0.$ On the other hand $t\rightarrow\frac{p_{n}(t)}{t}$ is continuous on
$]0,1]$ which implies that it is bounded on $[0,1]$ and then $(p_{n})_{n}\in
E_{u}$ and converges uniformly to $f$ . Consequently, $f$ $\in(E_{u})^{ru}.$
This means that $(E^{ru})_{u}\subset(E_{u})^{ru}$. Therefore $(E_{u}%
)^{ru}=(E^{ru})_{u}$. We obtain the desired result.
\end{proof}

\end{document}